\documentclass[11pt]{article}

\usepackage{latexsym}
\usepackage{amsfonts}
\usepackage{amsthm}
\usepackage{bezier}

\begin{document}

\newcommand{\CC}{\mathbb{C}}
\newcommand{\NN}{\mathbb{N}}
\newcommand{\QQ}{\mathbb{Q}}
\newcommand{\RR}{\mathbb{R}}
\newcommand{\ZZ}{\mathbb{Z}}

\newcommand{\iso}{\stackrel{\sim}{=}}
\newcommand{\tr}{{\rm tr}}

\newcommand{\ZA}{{{\rm ZA}}}
\newcommand{\ZB}{{{\rm ZB}}}
\newcommand{\WA}{{{\rm WA}}}
\newcommand{\WB}{{{\rm WB}}}
\newcommand{\TA}{{{\rm TA}}}
\newcommand{\TB}{{{\rm TB}}}

\newenvironment{bew}{Proof:}{\hfill$\Box$}

\newtheorem{satz}{Proposition}
\newtheorem{remark}[satz]{Remark}
\newtheorem{bsp}[satz]{Example}
\newtheorem{de}[satz]{Definition}
\newtheorem{lemma}[satz]{Lemma}
\newtheorem{kor}[satz]{Corollary}
\newtheorem{theo}[satz]{Theorem}
\newtheorem{satzdef}[satz]{Theorem and Definition}

\sloppy
\swapnumbers


\title{A braid group representation derived from handlebodies}
\author{Reinhard H\"aring-Oldenburg\\
{\small Mathematisches Institut, Bunsenstr. 3-5,
 37073 G\"ottingen, Germany}\\
{\small email: haering@cfgauss.uni-math.gwdg.de}}
\date{December 14, 1998}
\maketitle

\begin{abstract}
We define an action of Artin's braid group on a finite dimensional algebra.
\end{abstract}


\section{Introduction}

It seems that there is still interest in new representations of the braid group. 
One reason is the still open question if the braid group is linear.
Another is the search for finite dimensional braid algebras that generalize
the Birman-Wenzl algebra.  In this note we describe a representation that
is constructed in a new way. It uses a series of algebras associated to
handlebodies.

\section{The representation}

In this section we construct a braid group representation that 
is related to the action of the braid group on a free group.
So we start with the following well known definition and result:
\begin{satzdef}
The $n$-string braid group  $B_n$ is defined by generators
$X_1,\ldots,X_{n-1}$ and relations
\begin{eqnarray}
X_iX_jX_i&=&X_jX_iX_j \qquad |i-j|=1\\
X_iX_j&=&X_jX_i\qquad |i-j|>1
\end{eqnarray}
There is a monomorphism from the braid group into the
automorphism group of the free group in $n$
generators $\phi:B_n\rightarrow {\rm Aut}(F_n)$ given by \cite{burde}:
\begin{eqnarray}
\phi(X_i)(t_j)&:=& t_j\qquad j\neq,i,i+1\\
\phi(X_i)(t_{i+1})&:=& t_i\\
\phi(X_i)(t_i)&:=& t_it_{i+1}t_i^{-1}
\end{eqnarray}
\end{satzdef}

It is best to think of the free group $F_n$ as the fundamantal group of the 
handlebody of genus $n$, i.e. the product of a $n$-holed disc with the unit interval. 
The holes times the unit interval are called the axes of the handlebody and are
numbered by $1,\ldots,n$.
In recent work \cite{GHO} Temperley-Lieb algebras related 
to handlebodies have been studied. It is not necessary to give the
complete definitions used in that work. It suffices to consider the 
following algebras:
\begin{de}
Fix a natural number $n$ and parameters $c_{i,j}\in R,i,j\in\{1,\ldots,n\}$
in an integral domain $R$.
The algebra $T_n$ is defined to be generated by $f_{i,j},i,j\in\{1,\ldots,n\}$
with the set of relations:
\begin{equation}
f_{i,j}f_{k,l}=c_{j,k}f_{i,l}	  \label{rel1}
\end{equation}
\end{de}
It follows at once that $T_n$ is a free $R$-module of rank $n^2+1$.

In the handlebody setting we interpret $f_{i,j}$ as the $(1,1)$ tangle in the 
handlebody that is composed out of   a string that
enters axes $j$ and one that emerges from axes $i$. $c_{i,j}$ is then
the parameter that accounts for removing a string that connects 
two axes. Since Temperley-Lieb strings are not oriented one has
$c_{i,j}=c_{j,i}$. In the present paper we will find an even more 
restrictive condition on the parameters.

In accordance with the general philosophy of Temperley-Lieb 
algebras we interpret
\[Y_i:=\alpha_i f_{i,i}+\beta_i\qquad i\in\{1,\ldots,n\}\]
as the generator of the braid group of the handlebody that surrounds 
axes $i$. This means that $Y_i$ is to be considered as the image of $t_i$
under the Kauffman bracket map for the handlebody. 
The $\alpha_i,\beta_i$ are parameters in $R$ and we require $\alpha_i$ to be 
invertible.  This makes $Y_i$ invertible.
It is now clear, how to define the braid group action if we postulate, 
that the $c_{i,j}$ are invertible:

\begin{de} For each	  $i\in\{1,\ldots,n-1\}$ define operators 
$\phi(X_i)$ on $T_n$ by 
\begin{eqnarray}
\phi(X_i)(f_{j,j})&:=& f_{j,j}\qquad j\neq,i,i+1\\
\phi(X_i)(f_{i+1,i+1})&:=& f_{i,i}\\
\phi(X_i)(f_{i,i})&:=& \alpha_i^{-1}(Y_iY_{i+1}Y_i^{-1}-\beta_i)\\
\phi(X_i)(f_{i,j})&:=& \phi(X_i)(f_{i,i})\phi(X_i)(f_{j,j})c_{i,j}^{-1}\qquad i\neq j \\
\phi(X_i)(1)&:=&1
\end{eqnarray}
Denote by $B(i)$ the matrix of $\phi(X_i)$ with respect to the 
basis $1,f_{1,1},\ldots,f_{1,n},f_{2,1},\ldots,f_{n,n}$.
\end{de}

\begin{lemma}
The operators $\phi(X_i)$ are invertible. Inverses are given by
\begin{eqnarray}
\phi(X_i^{-1})(f_{j,j})&:=& f_{j,j}\qquad j\neq,i,i+1\\
\phi(X_i^{-1})(f_{i,i})&:=& f_{i+1,i+1}\\
\phi(X_i^{-1})(f_{i+1,i+1})&:=& \alpha_i^{-1}(Y_{i+1}^{-1}Y_{i}Y_{i+1}-\beta_i)\\
\phi(X_i^{-1})(f_{i,j})&:=& \phi(X_i^{-1})(f_{i,i})\phi(X_i^{-1})
(f_{j,j})c_{i,j}^{-1}\qquad i\neq j\\
\phi(X_i^{-1})(1)&:=&1
\end{eqnarray}
\end{lemma}
The proof is a simple calculation.
Note that the representation matrices $B(i)$ have off-diagonal
terms. 

It remains to establish that the $B(i)$ satisfy the relation of Artin's
braid group. To do this we will have to fix most of the parameters to
special values. We start by checking if the relation $f_{i,i}^2=c_{i,i}f_{i,i}$
is preserved under the action of $\phi(X_i)$. In the expansion of the result
the term proportional to $1$ is
\[\alpha_i^{-2}(\beta_i-\beta_{i+1})(c_{i,i}\alpha_i+\beta_i-\beta_{i+1})\]
We choose  the $\beta_i$ to be the same:
\begin{equation}
\beta_i:=\beta	\label{C1}
\end{equation}
Similarly, the coefficient of $f_{i,i}$ tells us that either $\alpha_i$
should vanish or that
\begin{equation} \label{C2}
  \alpha_i:= c_{i-1,i-1}\alpha_{i-1}/c_{i,i}\qquad  i>1
\end{equation}
Examining now the general relation (\ref{rel1})	yields
\begin{equation}
c_{i,j}:= c	 \label{C3}
\end{equation}

\begin{satz}
The map $\phi$ is a representation of the braid group if (\ref{C1}), 
(\ref{C2}), (\ref{C3}) hold.
\end{satz}
Again, this is a tedious calculation.

It is easy to see that the conditions on the parameters imply 
that the action of $\phi$ on $T_n$ preserves the subspace 
that is spanned by the $f_{i,j}$. The matrices given below are for this
subspace.

The representation depends in fact only on one parameter: 
\[ \mu=c^2\alpha_1/\beta   \]
From the calculation of examples we have the conjecture:
The characteristic polynomial ${\rm det}(B(i)-\lambda)$ is given by
\[
(\lambda-1)^{2+(n-1)^2}(\lambda+1+\mu)^{n-1}
(\lambda+(1+\mu)^{-1})^{n-1}
\]
The eigenvalues are thus $1,1+\mu,(1+\mu)^{-1}$ and the representation matrices have 
unit determinant. 
For the benefit of the reader, we give the matrices for the example $n=3$.
Experimental calculations show that the 3-string algebra they generate
is of dimension 19. Its centralizer algbra is of dimension $7$. Hence, 
if the algebra is semi-simple (as we hope) then it has $4$ two-dimensional and $3$ 
one-dimensional representations.

\[B(1)=\left(
\begin{array}{ccccccccc}
 -{\frac{{{\mu }^2}}{1 + \mu }} & {\frac{\mu }{1 + \mu }} & 0 & -
\mu  & 1 & 0 & 0 & 0 & 0 \\ \mu 
   & 0 & 0 & 1 + \mu  & 0 & 0 & 0 & 0 & 0 \\ 0 & 0 & {\frac{\mu }{1
 + \mu }} & 0 & 0 & 1 & 0 & 0
   & 0 \\ -{\frac{\mu }{1 + \mu }} & {\frac{1}{1 + \mu }} & 0 & 0 & 0 
& 0 & 0 & 0 & 0 \\ 1 & 0 & 0
   & 0 & 0 & 0 & 0 & 0 & 0 \\ 0 & 0 & {\frac{1}{1 + \mu }} & 0 & 0 
& 0 & 0 & 0 & 0 \\ 0 & 0 & 0
   & 0 & 0 & 0 & -\mu  & 1 & 0 \\ 0 & 0 & 0 & 0 & 0 & 0 & 1 + \mu
  & 0 & 0 \\ 0 & 0 & 0 & 0 & 
  0 & 0 & 0 & 0 & 1 \\  
\end{array}
\right)\]
\[B(2)=\left(
\begin{array}{ccccccccc}
 1 & 0 & 0 & 0 & 0 & 0 & 0 & 0 & 0 \\ 0 & -\mu  & 1 & 0 
& 0 & 0 & 0 & 0 & 0 \\ 0 & 
  1 + \mu  & 0 & 0 & 0 & 0 & 0 & 0 & 0 \\ 0 & 0 & 0 & {\frac{\mu }{
1 + \mu }} & 0 & 0 & 1 & 0 & 0
   \\ 0 & 0 & 0 & 0 & -{\frac{{{\mu }^2}}{1 + \mu }} & {\frac{\mu }{1 + 
\mu }} & 0 & -\mu  & 1 \\ 0 & 0
   & 0 & 0 & \mu  & 0 & 0 & 1 + \mu  & 0 \\ 0 & 0 & 0 & {\frac{1}{1
 + \mu }} & 0 & 0 & 0 & 0 & 0
   \\ 0 & 0 & 0 & 0 & -{\frac{\mu }{1 + \mu }} & {\frac{1}{1 + \mu }} 
& 0 & 0 & 0 \\ 0 & 0 & 0 & 0
   & 1 & 0 & 0 & 0 & 0 \\  
\end{array}
\right)\]

\section{Discussion}

We have found a braid group representation  from a simple 
geometric setting. This construction rises a number of open questions that should
be discusssed:
One should try to classify the irreducible subrepresentations. 

It is sensible to expect that the same method yields much more 
representations when applied
to other quotients of the handlebody braid group.

Another question: Is a Markov trace connected with this representation?

\end{document}